\documentclass[12pt]{article}
\usepackage{amsmath,latexsym,amssymb,graphicx,epsf,amsfonts,amsthm}
\usepackage{epsf}
\usepackage{color}
\numberwithin{equation}{section}

\begin{document}
\hoffset = -1truecm \voffset = -1truecm
\title{ A Semi-Blind Source Separation Method for Differential Optical Absorption Spectroscopy of Atmospheric Gas Mixtures }
\author{Y. Sun \thanks{Department of Mathematics,
University of California at Irvine, Irvine, CA 92697, USA.}, L.M.  Wingen
\thanks{Department of Chemistry,
University of California at Irvine, Irvine, CA 92697, USA.}, B.J. Finlayson-Pitts$^{\dagger}$, and J. Xin$^{*}$}
\date{}
\maketitle
\begin{abstract}
Differential optical absorption spectroscopy (DOAS) is a powerful tool for detecting and quantifying trace gases in atmospheric chemistry \cite{Platt_Stutz08}.  DOAS spectra consist of a linear combination of complex multi-peak multi-scale structures.  Most DOAS analysis routines in use today are based on least squares techniques, for example, the approach developed in the 1970s \cite{Noxon_75,Noxon_79,Per_76,Platt_79} uses polynomial fits to remove a slowly varying background (broad spectral structures in the data), and known reference
spectra to retrieve the identity and concentrations of reference gases \cite{Stutz_Platt96}.  An open problem \cite{Platt_Stutz08} is to
identify unknown gases in the fitting residuals for complex atmospheric mixtures.

In this work, we develop a novel three step semi-blind source separation method.
The first step uses a multi-resolution analysis to remove the slow-varying and fast-varying components
in the DOAS spectral data matrix $X$.  The second step decomposes the
preprocessed data $\hat{X}$ in the first step into a linear combination of the reference
spectra plus a remainder, or  $\hat{X} = A\,S + R$, where columns of matrix $A$ are known
reference spectra, and the matrix $S$ contains the unknown non-negative coefficients that are proportional to concentration.  The second step is realized by a convex minimization problem
$S = \mathrm{arg} \min \mathrm{norm}\,(\hat{X} - A\,S)$, where the norm is a hybrid $\ell_1/\ell_2$ norm (Huber estimator) that helps to maintain the non-negativity of $S$.  The third step performs
a blind independent component analysis of the remainder matrix $R$ to
extract remnant gas components.  We first illustrate the proposed method in
processing a set of DOAS experimental data by a satisfactory blind extraction
of an a-priori unknown trace gas (ozone) from the remainder matrix.  Numerical results also show that the method
can identify multiple trace gases from the residuals.
\end{abstract}

\thispagestyle{empty}
\newpage
\setcounter{equation}{0} \setcounter{page}{1}
\section{Introduction}
Trace gases play an important role in climate change and air quality of the Earth's atmosphere.  Spectroscopic techniques are widely used today for measurements of many trace species, and have evolved over the past century from the first use of the sun as a light source to identify atmospheric trace gases.  Many different light sources (e.g., infrared and UV-visible lamps, lasers, and natural sources such as the sun) are now conventionally used to identify light-absorbing species as well as determine their concentrations using Lambert-Beer's law,
\begin{equation}
\label{Lambert1}
I(\lambda) = I_0(\lambda)\cdot \exp(-\sigma(\lambda)\cdot \rho\cdot L)\;,
\end{equation}
where $I_0(\lambda)$ is the initial intensity of light, $I(\lambda)$ is its intensity after traveling through a sample of path length, $L$, with concentration, $\rho$.  Each species has its characteristic absorption cross section, $\sigma(\lambda)$, a measure of its ability to absorb light that varies with wavelength.  The use of (\ref{Lambert1}) is convenient for multi-component samples in laboratory spectrometers, but it is more difficult to determine the value of $I_0(\lambda)$ in the atmosphere over a large wavelength range.

A new method, differential optical absorption spectroscopy (DOAS), was introduced in the 1970s \cite{Noxon_75, Noxon_79, Per_76, Platt_79} to analyze atmospheric trace gas concentrations. DOAS analysis separates the trace gas absorptions, which typically vary quickly with wavelength, from features that vary slowly with wavelength, e.g., light scattering processes by molecules and aerosols.  Differential cross sections are then defined relative to this new broad background in place of the true $I_0(\lambda)$.  Several important trace gases were measured for the first time with DOAS, e.g., HONO, $\mathrm{NO}_3$, BrO, ClO in the troposphere, and OClO and BrO in the stratosphere.  A large number of other molecules absorb in the UV and the visible wavelength region and most aromatic hydrocarbons can also be detected.  An advantage of DOAS is the ability to measure absolute trace gas concentrations {\it in situ}.  DOAS is therefore especially useful for measuring highly reactive species such as the free radicals OH, $\mathrm{NO}_3$, or BrO, and it provides a powerful tool for studying emissions, transformation and transport of chemicals throughout the troposphere and stratosphere. It can also help to understand the influence of atmospheric chemistry on climate and air quality.  A detailed description of DOAS can be found elsewhere \cite{Platt_Stutz08}.

In general, DOAS spectra contain overlapping absorption structures which consist of complex multiple scales and peaks.  They must be separated by the analysis routine to retrieve the concentrations of the trace gases.  Least squares techniques are most often used for analysis of DOAS spectra, with the use of high pass filters to fit or separate out the slowly varying components.  For example, the approach described in \cite{Platt_Stutz08} applies a polynomial fit to remove the broad (slow-varying) spectral features, and known reference spectra to retrieve the concentrations of reference gases.  However, the existing DOAS approaches have two limitations: 1) the condition of least squares ({\em that errors are normally distributed} )
is often violated.  This suggests that a different norm other than $\ell_2$ (least squares) norm should be used;
2) the fitting residuals for atmospheric samples are in most cases not pure noise due to
imperfect references, atmospheric turbulence, instrument effects, and unknown trace gases.  Among other interesting problems, the identification of gas structures in the fitting residuals is of great importance.
The method in this paper has been developed to address these issues, in the hope of
providing a tool for atmospheric chemists to analyze the residuals for possible hidden trace gases.
The method is designed to deal with the following three major challenges.
First, DOAS spectra are complex multi-scale multi-peak structural data containing slow-varying features, structured signals due to the trace gases, and instrumental noise.  Hence a multi-resolution analysis tool
is needed for scale decomposition.  Second, the identification of gases from the residuals is actually a problem of
blind source separation (BSS) as both the trace gases (including their numbers)
and mixing process are not known.  A major problem is to find a working assumption on the
source (hidden trace gas) signals and effective BSS algorithms.
Third, the new objective function for data fitting should not only overcome
the limitations of least squares fitting, but also help to maintain the non-negativity.
To tackle these problems, we have made an initial attempt of developing a semi-blind approach
which contains three steps.  The first step uses  multi-resolution analysis to
remove the very slow (e.g. scattering) and very fast components (noise) in the DOAS spectral data matrix $X$.
The second step decomposes the preprocessed data $\hat{X}$ in the first step into a linear combination
of the reference spectra plus a remainder, or  $\hat{X} = A\,S + R$, where columns of matrix $A$ are known
reference spectra, the matrix $S$ contains the unknown non-negative coefficients.
The second step is carried out by solving a convex minimization problem
$S = \mathrm{arg }\min \mathrm{norm}\,(\hat{X} - A\,S)$, where the norm is a hybrid $\ell_1/\ell_2$ norm that
helps to maintain the non-negativity of $S$.  The third step performs
a blind independent component analysis of the remainder matrix $R$ to
extract remnant gas components.  Our method can be useful for separating unknown sources
from the residuals after any known reference spectra have been first deployed to fit the data.

The paper is organized as follows.  In section 2, we review the
essentials of DOAS and the existing approach, then introduce our method.  In section 3,
we illustrate the proposed method in processing a set of DOAS experimental data, and
show satisfactory numerical results.  Concluding remarks are in section 4.
\section{DOAS and Signal Processing Methods}
\subsection{DOAS and Fitting Methods}
A typical experimental setup for a DOAS instrument consists of a continuous light source, e.g., a Xe-arc lamp,
a light-absorbing sample (the atmosphere or gases in a chamber), a grating spectrometer, and an optical detector.  It is also possible to use the light from the sun or moon, or scattered sun light as light sources \cite{Noxon_75, Noxon_79}. The typical length of the light path in the atmosphere ranges from several hundred meters to many kilometers and $<$ 100 m in laboratory DOAS experiments.  The light of intensity $I_0(\lambda)$ passes through the sample, is typically  dispersed by a grating spectrometer and is measured by a detector.  During its way through the sample the light undergoes extinction due to absorption processes by trace gases and scattering by air molecules and aerosol particles.  In the atmosphere, the intensity $I(\lambda)$ at the end of the light path is given by Lambert-Beer's law,
\begin{equation}
\label{Lambert}
I(\lambda) = I_0(\lambda)\exp \Biggl [ -\int^{L}_{0}  \sum_{j} \sigma^{\mathrm{ABS}}_j(\lambda) \times \rho_j(l) + \varepsilon_{\mathrm{R}}(\lambda,l) + \varepsilon_{\mathrm{M}}(\lambda,l)\; \mathrm{d}l \Biggr] + N(\lambda)\;,
\end{equation}
where $\sigma^{\mathrm{ABS}}_j$ is the absorption cross section of a trace gas $j$,  $\rho_j$ is its number density.  $L$ is the length of the light path.  The Rayleigh extinction by gases and Mie extinction by aerosols are described here by $\varepsilon_{\mathrm{R}}$ and $\varepsilon_{\mathrm{M}}$.  $N(\lambda)$ is the measurement noise.  The basic idea of DOAS is the separation of the cross section $\sigma^{\mathrm{ABS}}_j = \sigma^{\mathrm{B}}_j + \sigma'_j$ in which  $\sigma^{\mathrm{B}}_j$ represents broad spectral features and the differential cross section $\sigma'_j$ represents narrow spectral structures that are of interest for identification and quantification of the trace gases.  If one considers only $\sigma'_j$, interferences with Rayleigh and Mie extinction are avoided.
The mathematical description of this process is a convolution of $I(\lambda)$ with the instrument function $H$ of the spectrometer,
\[ I^*(\lambda)  =  I(\lambda)*H =  \int I(\lambda - \lambda')\; H(\lambda')\mathrm{d}\lambda' = \]
{\allowdisplaybreaks
\begin{eqnarray*}
  \int^{\Delta \lambda}_{-\Delta \lambda} I_0(\lambda-\lambda')\exp \Biggl [ -\int^{L}_{0}  \sum_{j} \sigma^{\mathrm{ABS}}_j(\lambda-\lambda')\; \rho_j(l) + (\varepsilon_{\mathrm{R}}+ \varepsilon_{\mathrm{M}})(\lambda-\lambda',l)\, \mathrm{d}l \Biggr] \; H(\lambda') \mathrm{d}\lambda'\\
   =     \int^{\Delta \lambda}_{-\Delta \lambda} I'_0(\lambda-\lambda')\exp \Biggl [ -\int^{L}_{0}  \sum_{j} \sigma'_j(\lambda-\lambda') \times \rho_j(l) \mathrm{d}l \Biggr] \; H(\lambda') \mathrm{d}\lambda'\;,
\end{eqnarray*}}
where $$\displaystyle I'_0(\lambda-\lambda') = I_0(\lambda-\lambda')\exp \Biggl [ -\int^{L}_{0}  \sum_{j} \sigma^{\mathrm{B}}_j(\lambda-\lambda') \;\rho_j(l) + (\varepsilon_{\mathrm{R}}+ \varepsilon_{\mathrm{M}})(\lambda-\lambda',l)\;\mathrm{d}l \Biggr] $$
describes the broad spectral structures due to the characteristics of the light source $I_0$, the Rayleigh and Mie's extinction, and the broad absorption by trace gases.  $I'_0$ is a slow-varying function of wavelength, so $I^*(\lambda)$ can be approximated by
{\allowdisplaybreaks
\begin{eqnarray*}
I^*(\lambda)& = & I'_0(\lambda) \int^{\Delta \lambda}_{-\Delta \lambda}\exp \Biggl [ -\int^{L}_{0}  \sum_{j} \sigma'_j(\lambda-\lambda') \times \rho_j(l) \mathrm{d}l \Biggr] \times H(\lambda') \mathrm{d}\lambda' \\
& =&  I'_0(\lambda) \exp \Biggl [  S_j \times \int^{\Delta \lambda}_{-\Delta \lambda} \sum_{j} -\sigma'_j(\lambda-\lambda') \; \mathrm{d}l \Biggr] \times H(\lambda') \mathrm{d}\lambda'\\
& = & I'_0(\lambda) \exp \biggl [\sum_{j}S_j \times   A_j(\lambda) \biggr ]
\end{eqnarray*}}
where $\displaystyle S_j = \int^L_0 \rho_j(l)\;\mathrm{d}l $, and
$A_j(\lambda) = -\sigma'_j(\lambda)* H(\lambda)$ denote
the narrow absorption structures of the trace gases measured with the same instrument.
Suppose that there are $m$ known trace gases in the data.  The logarithm of the above equation becomes
\begin{equation}
\label{log}
x(\lambda) = \sum^{m}_{j = 1} S_j \times A_j(\lambda) + B'(\lambda) + N'(\lambda)\;,
\end{equation}
where $\displaystyle x(\lambda) = \ln I^*(\lambda)$ and $B'(\lambda) = \ln I'_0(\lambda)$ represents
the broad spectral features.  The noise $N'(\lambda) = \ln N(\lambda)$.  In the experiment, the wavelength range is mapped onto $p$ discrete pixels of the detector.  The sampled data points form $p$-dimensional column vectors of a data matrix $X$.
Suppose there are $n$ measurements, then
 $X = [x_1, x_2, \dots, x_n] \in \mathbb{R}^{p\times n}$ whose column vectors are
recorded DOAS data points.  Equation (\ref{log}) in  matrix form is:
\begin{equation}
\label{modelEQ}
X = AS + B + N\;,
\end{equation}
where the columns of matrix $A$ correspond to the reference spectra of the known trace gases;
the matrix $S$ contains non-negative coefficients;
the matrix $B$ includes the slow-varying components, and the matrix $N$ contains the noise components.
Most DOAS approaches use the least squares methods to calculate $S$ due to its computational simplicity.
To deal with the problem that the data contain both broad and narrow spectral features,
a high pass filter is needed to remove the broad spectral features.  It is common to use polynomials to model and filter out the slowly varying parts
from the narrow trace gas absorption.  Equation (\ref{modelEQ}) is written as follows in \cite{Platt_Stutz08}.
\begin{equation}
\label{poly}
X = AS + P + N\;,
\end{equation}
where the polynomial $P$ models the broad spectral structures.
Given the order of polynomial and the known reference spectra, (\ref{poly}) can be solved
with a least squares method.
The polynomial fitting however
has the following drawbacks: (1) the order of the polynomial is determined empirically
and different orders might be used for different data; 2)
the non-negativity of the concentration is not guaranteed during the fitting.
An open problem after the fitting is
how to identify and extract trace gases from the fitting residuals besides the noise.
To address these issues, we propose a three step method in the next section.
\subsection{Proposed Semi-Blind Source Separation Method}
DOAS data can cover a range of scales and contain high frequency ($<$1 nm) artifact structures, for example due to pixel-to-pixel variability in the detector, while the reference spectra of the trace gases contain fewer peaks and peak widths on the order of several nm.  Hence it is helpful to remove the fast varying artifacts from the spectra data by multi-resolution analysis.  In addition, the broad features (slow-varying parts) in the data need to be eliminated
in order to fit the reference spectra of the known trace gases.
We propose to use the empirical mode decomposition (EMD) to extract these components.
The detailed description of EMD can be found in \cite{Huang_98}.  The EMD method does not assume anything about the data, contrary to Fourier methods where data is assumed linear and stationary.  EMD handles also non-stationary and nonlinear data.
 
\subsubsection{Multi-Resolution Analysis}
The concept of EMD has been developed
rapidly in many areas of science and
engineering since Huang et al. \cite{Huang_98} invented EMD.
Its key feature is to decompose a signal
into so-called intrinsic mode function (IMF).  The essential step extracting an IMF is to identify
an oscillation embedded in a signal from local
time scale.  Considering a signal $s(t)$ between two consecutive local
extrema (say, two minima at times $t_1$ and $t_2$),  we can heuristically
define a (local) high frequency part $\{d(t), t_1\leq t_2 \}$, where
$d(t)$ corresponds to the oscillation terminating at the two minima and
passing through the maximum in between.  For the picture to be complete,
we also identify the corresponding (local) low-frequency part, or local {\sl trend},
 $m(t)$ so that we have $s(t) = m(t) + d(t) $ for $t_1\leq t\leq t_2$.
Assuming that this is done in some proper way for all the oscillations in the entire signal,
we get an intrinsic mode function as well as
a residual consisting of all local trends. The procedure can then be repeated on the residual, and constitutive components of a signal can be iteratively extracted.

\begin{figure}
\includegraphics[height=7.5cm,width=16cm]{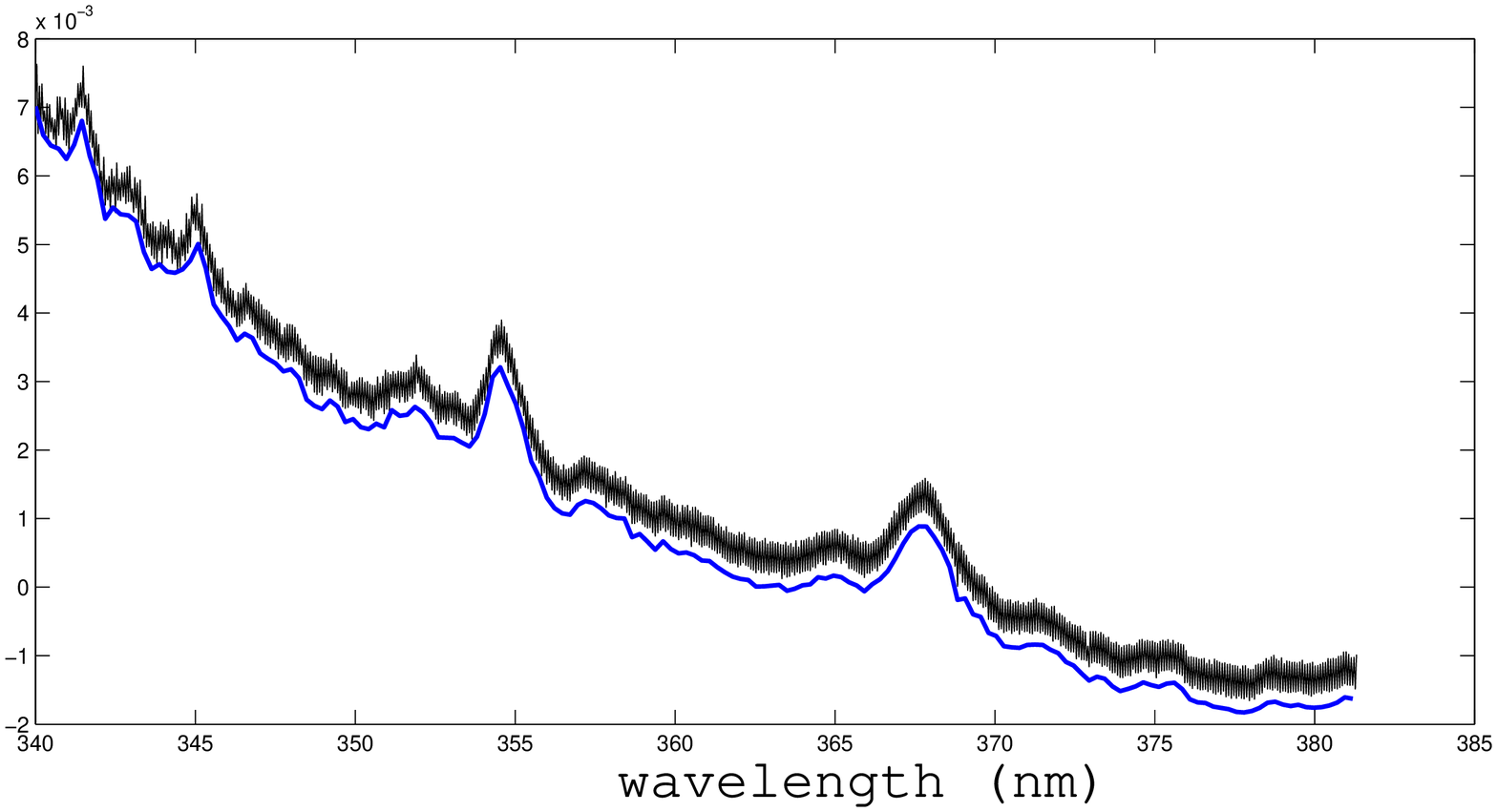}
\includegraphics[height=7.5cm,width=16cm]{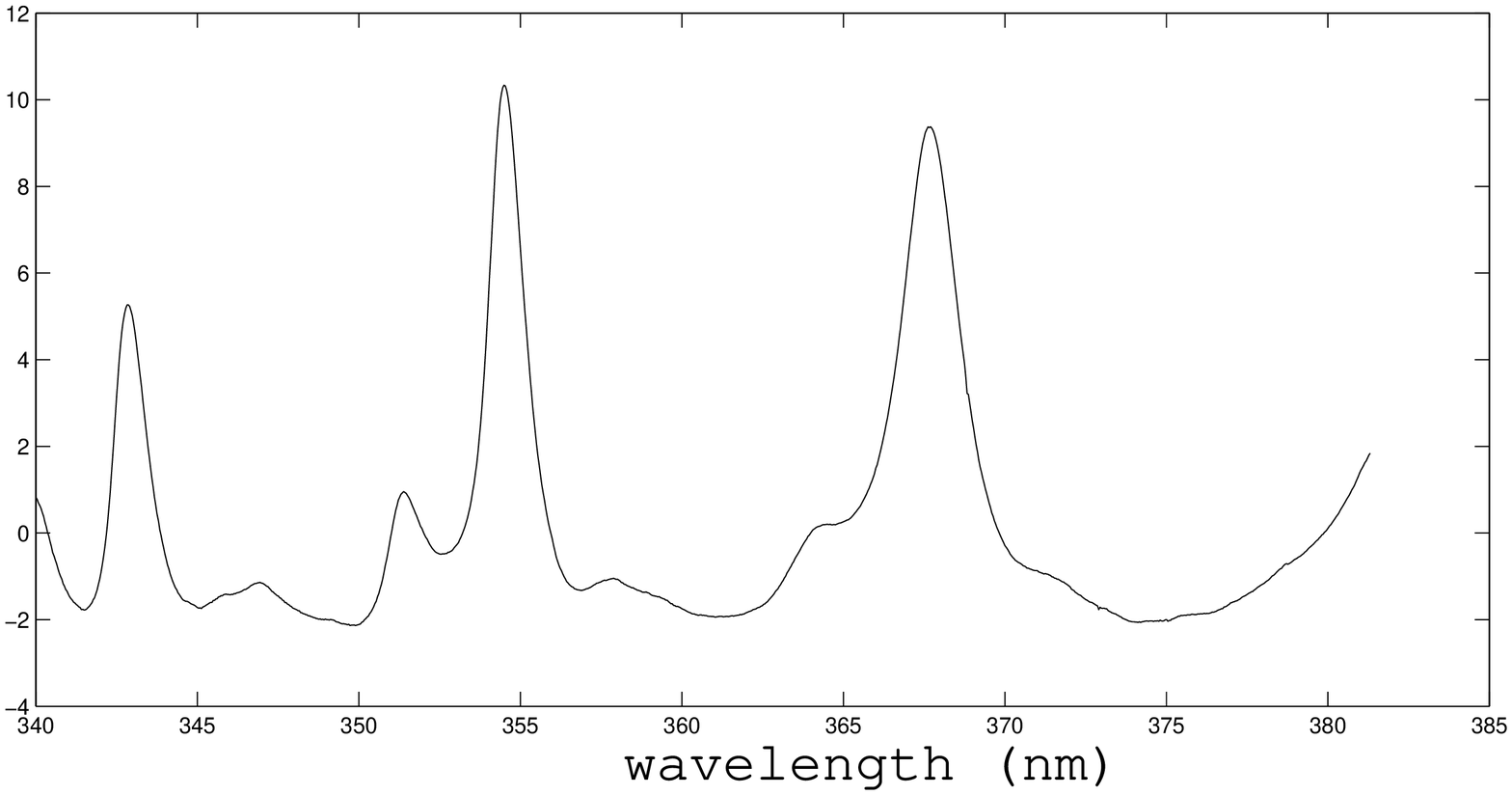}
\caption{Top panel is a mixed DOAS spectrum (black) from the experiment and its smoothed counterpart (blue).
Bottom panel is the spectral absorption reference of trace gas HONO.}
 \label{DOSAdata1}
\end{figure}

\begin{figure}
\includegraphics[height=8.5cm,width=17cm]{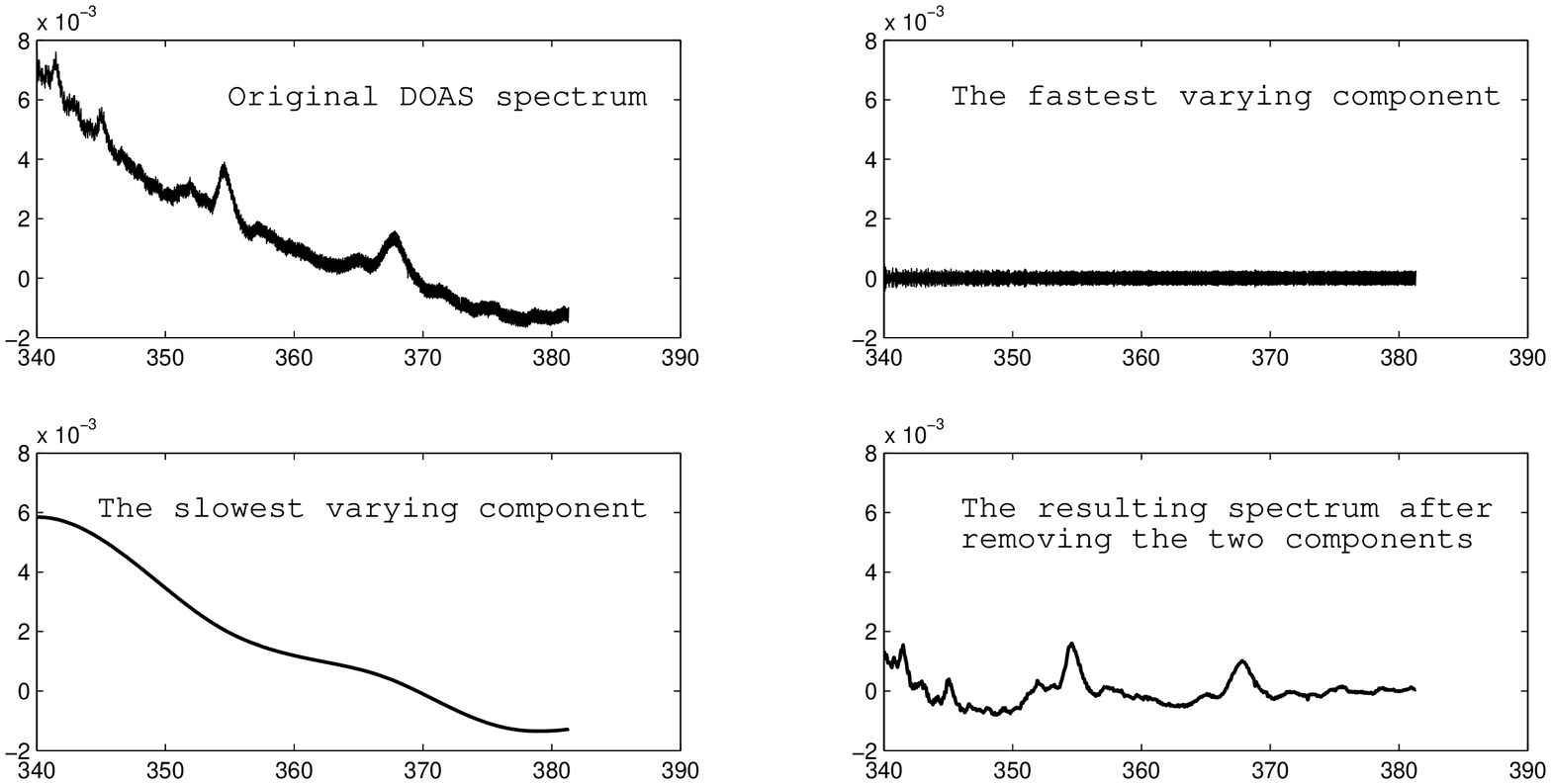}
\caption{The preprocessed DOAS data after removing the slowest and fastest varying components.  The x-axis of the plots is wavelength (nm). }
 \label{result1}
\end{figure}

The EMD method decomposes the DOAS data into a finite number of components of different frequencies.
The advantage of EMD is that it is completely data-driven (no need to specify a parameter such as
the order of a fitting polynomial), fast and
automatic.  Fig. \ref{DOSAdata1} shows a typical DOAS spectrum  of a trace gas mixture containing HONO, whose reference spectrum is also shown in the bottom panel.  Fig. \ref{result1} shows the fast and
slow components extracted from the DOAS data in Fig. \ref{DOSAdata1}.
The EMD preprocessed (high-passed) data $\hat{X}$ satisfies the following model
\begin{equation}
\label{model}
\hat{X} =  A\,S + R,
\end{equation}
where the columns of matrix $A$ are the reference spectra of the known trace gases,
and those of $S$ matrix contains their concentrations, and $R$ is the fitting residual which might contain the instrument noise, hidden trace gas structures, etc. For the estimation of the concentration matrix $S$, we minimize the following constrained objective function:
\begin{equation}
\label{convexModel}
\min_{S} \mathrm{norm}\,(\hat{X}-A\,S),\;\; \mathrm{s.t.}\;\; S\geq 0,
\end{equation}
for a proper choice of the norm.
\medskip

\subsubsection{Huber Estimator and Robust Data Fitting}
\begin{figure}
\includegraphics[height=7.5cm,width=15cm]{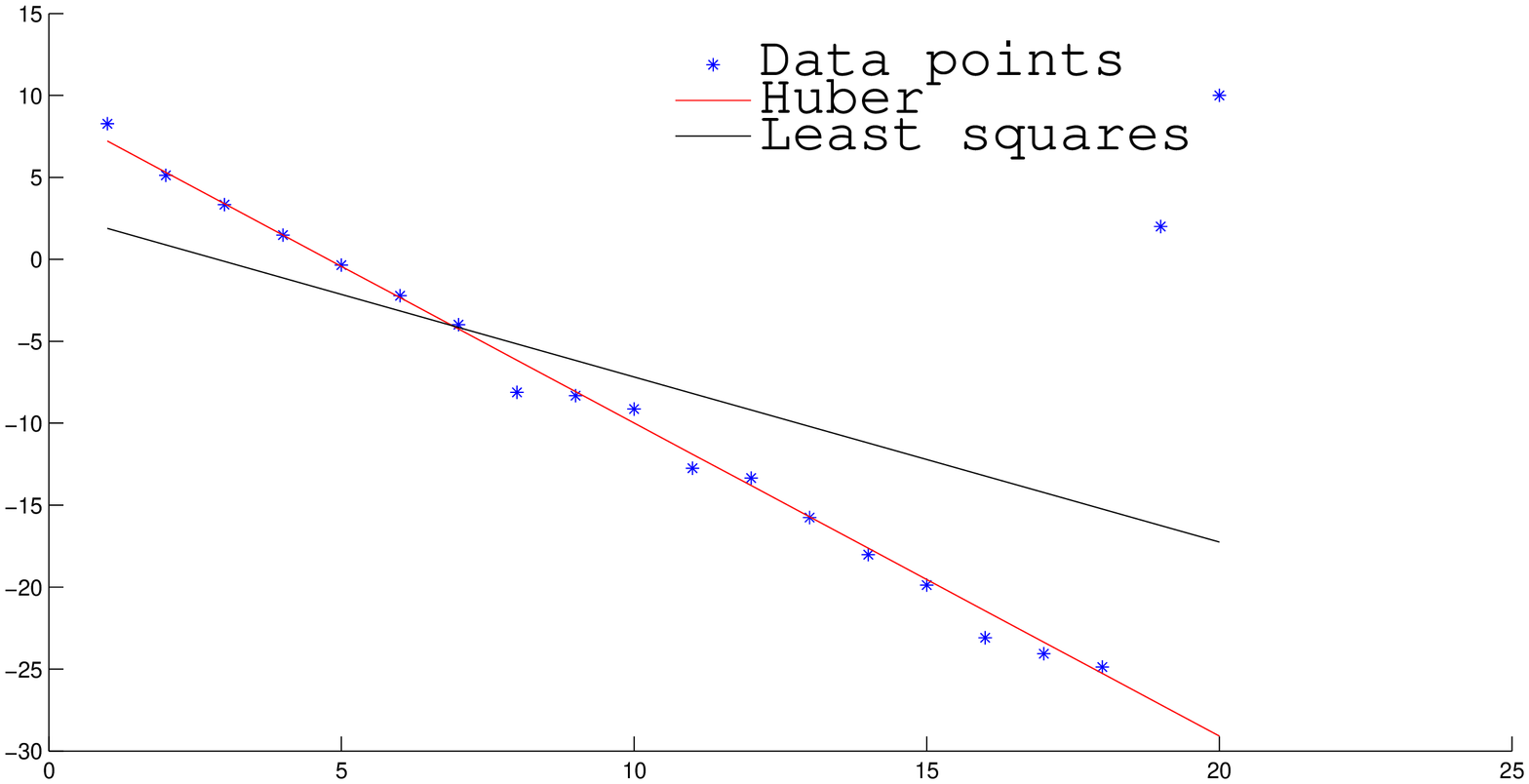}
\includegraphics[height=7.5cm,width=17cm]{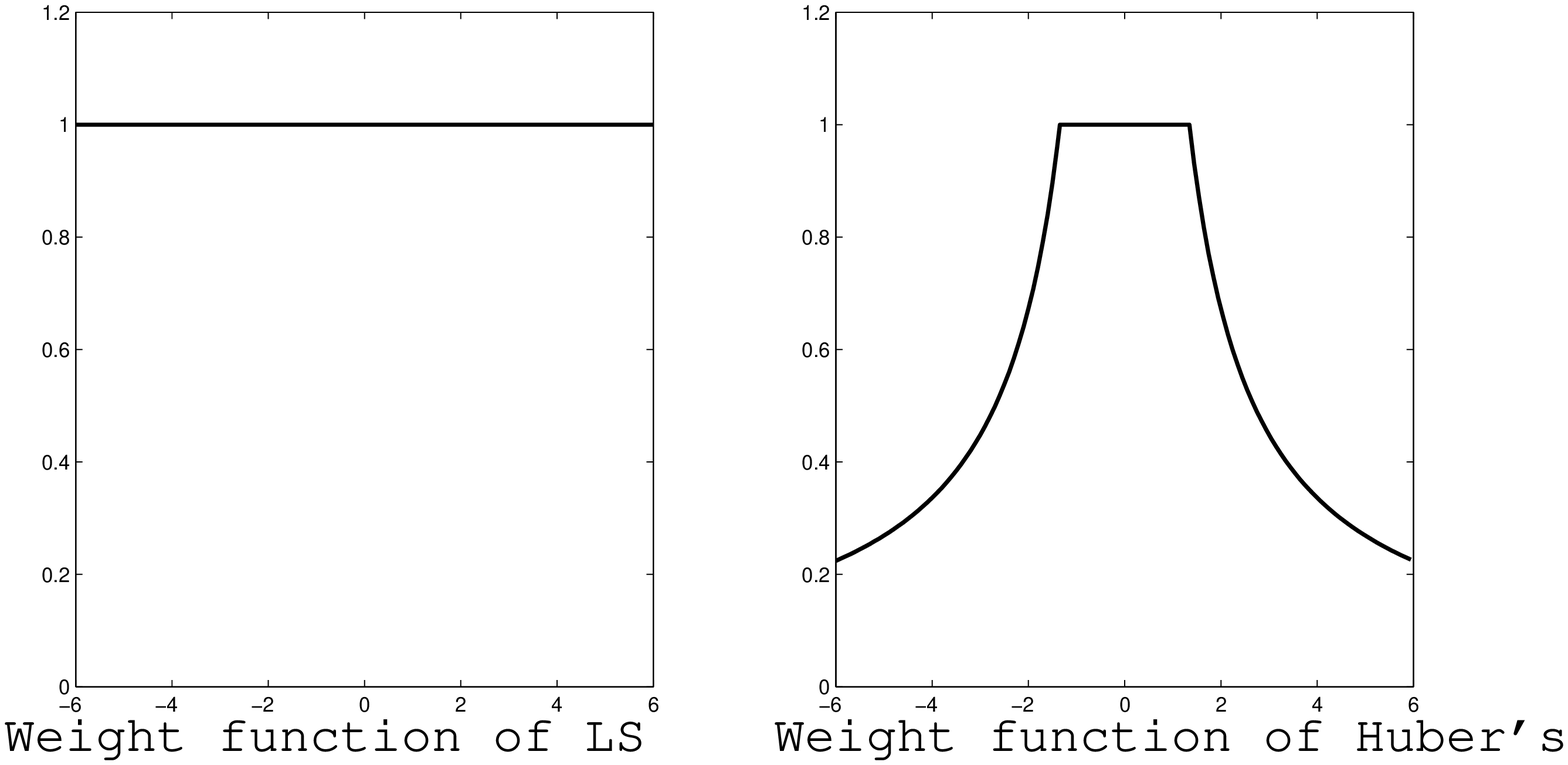}
\caption{Comparison of Huber's estimator and least squares (their weight functions are shown in the bottom plots). The data points are generated by $y = -2x + 10$ plus noise.  Huber: $y = -1.9794x + 9.9318$;
 Least squares: $y = -1.0504x + 3.5819$.}
 \label{huber_example}
\end{figure}

There are many kinds of norms available, e.g., $\ell_2$ (least squares), $\ell_1$ (least absolute deviations).
The regular least squares method (ignoring the non-negative constraint on $S$) is the conventional choice,
if the unknown noise $N$ is assumed to be Gaussian.
However, it is rather sensitive to the outliers in the data, even one outlier can drastically change the
estimation.  The least absolute deviations ($\ell_1$ norm) is more robust
to outliers in the data, however it is less effective if the peaks in
$N$ are not isolated (or sparse).  We find that a hybrid $\ell_2$ and
$\ell_1$ norm ($\ell_2$ on small peaks and $\ell_1$ on large peaks),
or a Huber estimator \cite{Huber_09}, is able to resist the influence of outliers,
and maintain non-negativity of $S$ for our data fitting task.
The Huber norm is both regular and convex. The corresponding nonlinear function $H=H(x)$ is a
parabola ($\ell_2$) in the vicinity of zero,
and increases linearly ($\ell_1$) at $|x| > k$ for any positive constant $k$. More precisely,
\begin{equation}
\label{huber}
H(x) = \left\{
\begin{array}{ll}
 \frac{1}{2}x^2,\;\;\;\;\;  |x| \leq k\\
 k|x| - \frac{1}{2}k^2,\;\;\;\;\; |x| \geq k.
  \end{array}
 \right.
\end{equation}
The non-negativity of $S$ under the Huber norm
indicates that the choice fits the empirical distribution of the noise $N$ arising from detection and
photon statistics \cite{Platt_Stutz08}.
Fig. \ref{huber_example} is an example showing the superiority of Huber's estimator
over least squares: it is resistant to outliers in the data, while the
least squares result deviated significantly from the exact line due to the two outliers.
Least squares assigns equal weighting to each observation;
the weights for the Huber estimator decline when $|x| > k$ (see the weight functions
in Fig. \ref{huber_example}).
The value $k$ for the Huber estimator is called a tuning constant;
smaller values of $k$ produce more resistance to outliers, but at the expense of
lower efficiency when the errors are normally distributed.
The tuning constant is generally selected to yield high efficiency in the normal case;
in particular, $k = 1.345\sigma$ (where $\sigma$ is the standard deviation of the errors).
The minimization of Huber's objective can be achieved by the method of iterative re-weighted least squares.
For the DOAS data in our numerical experiments, we found that
the Huber's estimator produced non-negative solutions $S$ without explicit enforcement
of non-negativity constraint.  The theoretical underpinning is under further study.
In the residual of Huber estimation, there might be spectral structures (one or many) of
trace gases buried in noise, or just random noise.  In either case,
we decompose the residuals in a blind fashion due to the lack of the knowledge of the hidden trace gases.
The source signal assumption required for the decomposition is that the spectra of different trace gases
are statistically independent (orthogonal). This appears to be a reasonable working assumption for
many trace gases.  Independent component analysis (ICA) can now be readily applied.
\subsubsection{Independent Component Analysis}
ICA is a useful and generic tool for solving blind source separation problems (BSS),
which arise when one attempts to recover source signals from their mixtures
without knowing the mixing process \cite{Cic_05,Comon_10}.
ICA finds the independent components in the mixtures
by maximizing the statistical independence (minimizing mutual information) of the estimated components.
Mathematically, given the mixture matrix $R\in \mathbb{{R}}^{p\times n}$ and
the number of independent source components $d$,
ICA finds a full rank matrix $W \in \mathbb{R}^{d\times n}$ such that the output matrix
$U \in \mathbb{{R}}^{d\times p}$ given by
\begin{equation}
\label{ICA}
U = W\;R'
\end{equation}
contains columns (recovered source signals) as independent from each other as possible.
Here $n$ is the number of residuals from the data fitting,
$p$ is the number of wavelength pixels.
The columns of $U$ correspond to the decomposed independent source signals.
We may choose one of many ways to approximate independence,
and this choice governs the form of the ICA algorithm.
The two broadest definitions of independence for ICA are:
(1) minimization of mutual information; (2) maximization of non-Gaussianity.
The non-Gaussianity family of ICA algorithms
use kurtosis and negentropy.  The minimization of mutual information family
of ICA algorithms use the Kullback-Leibler divergence and maximum-entropy, however,
the knowledge of source signal probability distribution function (PDF) is needed.
Algorithms for ICA include infomax \cite{Infomax}, FastICA  \cite{Fastica}, and JADE  \cite{Jade}.
We opt for JADE because JADE is based on cumulants (2nd and 4th order statistics)
and the approximate joint diagonalization of cumulant matrices (hence does not rely on PDF information of
source signals). For moderate number of sources,
it is more direct and stable
than iterative methods such as
infomax \cite{Infomax} and FastICA  \cite{Fastica}.
It was recently found \cite{LXQ} that the infomax method \cite{Infomax}  may even diverge and that
it only converges in a weak sense under proper rescaling and soft dynamic control of the iterations.
The most attractive aspect of JADE is that it does not require parameter tuning (e.g. choosing the
learning parameter in the iterative methods). In general, ICA algorithms cannot
identify the actual number of source signals, so this number needs to be found by other means, for example by
human evaluation of the end results.
In our decomposition of Huber residuals, we tested a range for this number, and
pinpointed the one with the most reliable and meaningful outcomes when calibrated with the knowledge of
the existing trace gas spectral properties.

\section{Experiments and Computational results}
\subsection{Experimental Setup}
\begin{figure}
\begin{center}
\includegraphics[height=10cm,width=8cm]{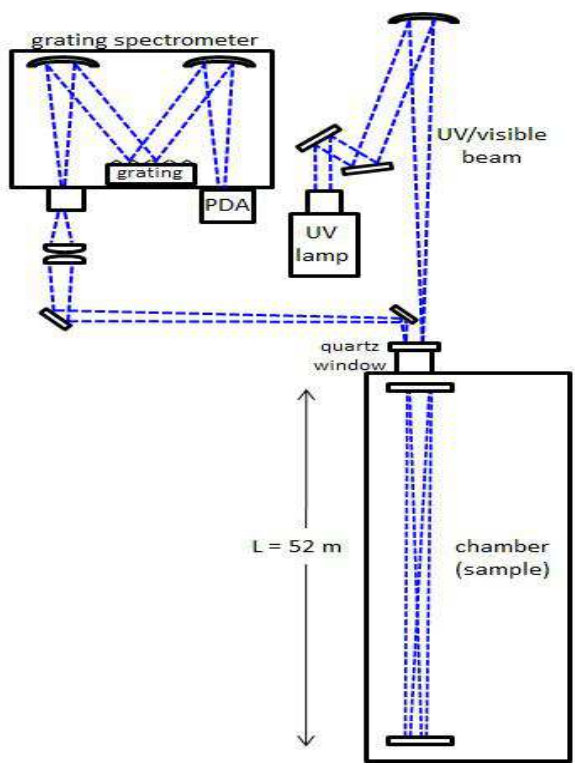}
\caption{Schematic of chamber, instrumentation, and optical setup used to make DOAS measurements of gaseous mixtures.}
 \label{doas_setup}
 \end{center}
\end{figure}

Spectra of chemical mixtures were collected using an environmental chamber \cite{Deh_99} for which DOAS is one of the analytical techniques used to measure species during experiments.  Fig. \ref{doas_setup} shows a simplified schematic of the chamber and optical arrangement for DOAS.  The chamber is 561 $L$ in volume and can be evacuated to a pressure of $ \sim 10^{-2}$ Torr for collection of true $I_0(\lambda)$ spectra.  Spectra can also be collected before and after addition of ultrapure air and each gaseous analyte of interest through various ports.
\medskip

The DOAS instrumentation consists of a high pressure Xe arc lamp (Oriel, Model 6263) as the UV-visible light source.  The light beam enters the chamber through a quartz window and undergoes multiple reflections using White cell mirrors \cite{White_42} through the gas mixture in the chamber.  The multiple reflections increase the path length of the light beam through the sample to a total path length of $L$ = 52 m.  The light beam exits the chamber through the quartz window and is focused on the entrance slit of a monochromator (Jobin Yvon-Spex, Model HR460) with a diode array detector (Princeton Instruments, model PDA-1024 ST121).  The grating (1200 grooves $\mathrm{mm}^{-1}$ blazed at 330 nm) gives a dispersion of $\sim 0.043$ \,nm $/\mathrm{pixel}$ and the detector has 1024 channels giving each spectrum a total wavelength range of $\sim$ 44 nm.  Spectra can be collected in different wavelength ranges by moving the grating motor.  Changes in grating position as well as temperature lead to changes in dispersion of the light beam on the detector.  This is taken into account in the least squares analysis by allowing for shifting or linear compression/expansion in one or more reference spectra along the wavelength axis to obtain the best fit.  The use of such techniques is standard and user controlled to correlate wavelengths with channels of the detector.  Absolute dispersion and wavelength were calibrated using a mercury lamp spectrum that was recorded daily and  at the beginning of each experiment.
\medskip

The analytes added to the chamber were $\mathrm{NO}_2$ and $\mathrm{O}_3$ at a total pressure of $ \sim 1$ atm at room temperature in dry ultrapure air (Scott-Marrin, Riverside, CA).  The wavelength range typically used to measure $\mathrm{NO}_2$ by DOAS is 340 - 380 nm.  Although the air was dry (relative humidity $< $ 0.8\%), even small amounts of water react with $\mathrm{NO}_2$ to form HONO \cite{BJFP_LMW_03}.  As a result HONO is almost always present in detectable quantities with $\mathrm{NO}_2$.  HONO is also typically measured using the 340 - 380 nm wavelength range, thus the mixture of HONO and $\mathrm{NO}_2$ was used as a convenient test case for the new DOAS analysis technique.  The addition of $\mathrm{O}_3$ leads to formation of $\mathrm{NO}_3$ radicals ($\mathrm{NO}_2$ + $\mathrm{O}_3$ $\rightarrow$ $\mathrm{NO}_2$ + $\mathrm{O}_2$).  Analysis for $\mathrm{O}_3$ is typically carried out in a different wavelength range, 290 - 330 nm.  It should be noted that a wavelength range for analysis is usually that in which the cross sections are highest for that analyte in order to optimize the detection limits.  $\mathrm{O}_3$ continues to absorb at wavelengths $>$ 330 nm, albeit with absorption cross sections that are lower by a factor of 100 or more \cite{Voigt_01} compared to those at shorter wavelengths.  Another test of the technique was to determine if it could identify the presence of this third component, $\mathrm{O}_3$, in the 340 - 380 nm range where its detection is not optimal.  $\mathrm{NO}_3$ analysis was carried out in a different range, 600 - 640 and 640 - 680 nm, and is not discussed here.
\medskip

In addition to the new DOAS analysis technique introduced in this work, the typical linear least squares analysis was carried out on HONO and $\mathrm{NO}_2$ using MFC \cite{Gomer_95} for which reference spectra are needed.  A reference spectrum for $\mathrm{NO}_2$ was generated by adding a known quantity of $\mathrm{NO}_2$ to the chamber and collecting DOAS spectra with the instrumentation described above.  Pure samples of HONO are difficult to generate without the presence of $\mathrm{NO}_2$, thus HONO reference spectra were generated from published cross sections \cite{Bongartz_91, Bongartz_94} which were convoluted to the dispersion and resolution of our spectrometer.  Reference spectra for $\mathrm{O}_3$ were generated from published cross sections \cite{Voigt_01} also converted to the dispersion and resolution of the spectrometer.  

\medskip
Chemicals used in these experiments are as follows:  Gaseous $\mathrm{NO}_2$ was synthesized by reaction of gaseous NO (Matheson, 99\%) which was first passed through a cold trap at 195 K to remove impurities such as $\mathrm{HNO}_3$, with an excess of $\mathrm{O}_2$ (Oxygen Services Co., 99.993\%).  The mixture was allowed to react for ~2 hrs. and then purified by condensing the $\mathrm{NO}_2$ at 195 K to pump away excess $\mathrm{O}_2$.  Gaseous $\mathrm{O}_3$ was generated as a mixture in $\mathrm{O}_2$ using a commercial ozonizer (Polymetrics, Model T-816).

\subsection{Computational Results}
\begin{figure}
\includegraphics[height=6cm,width=8cm]{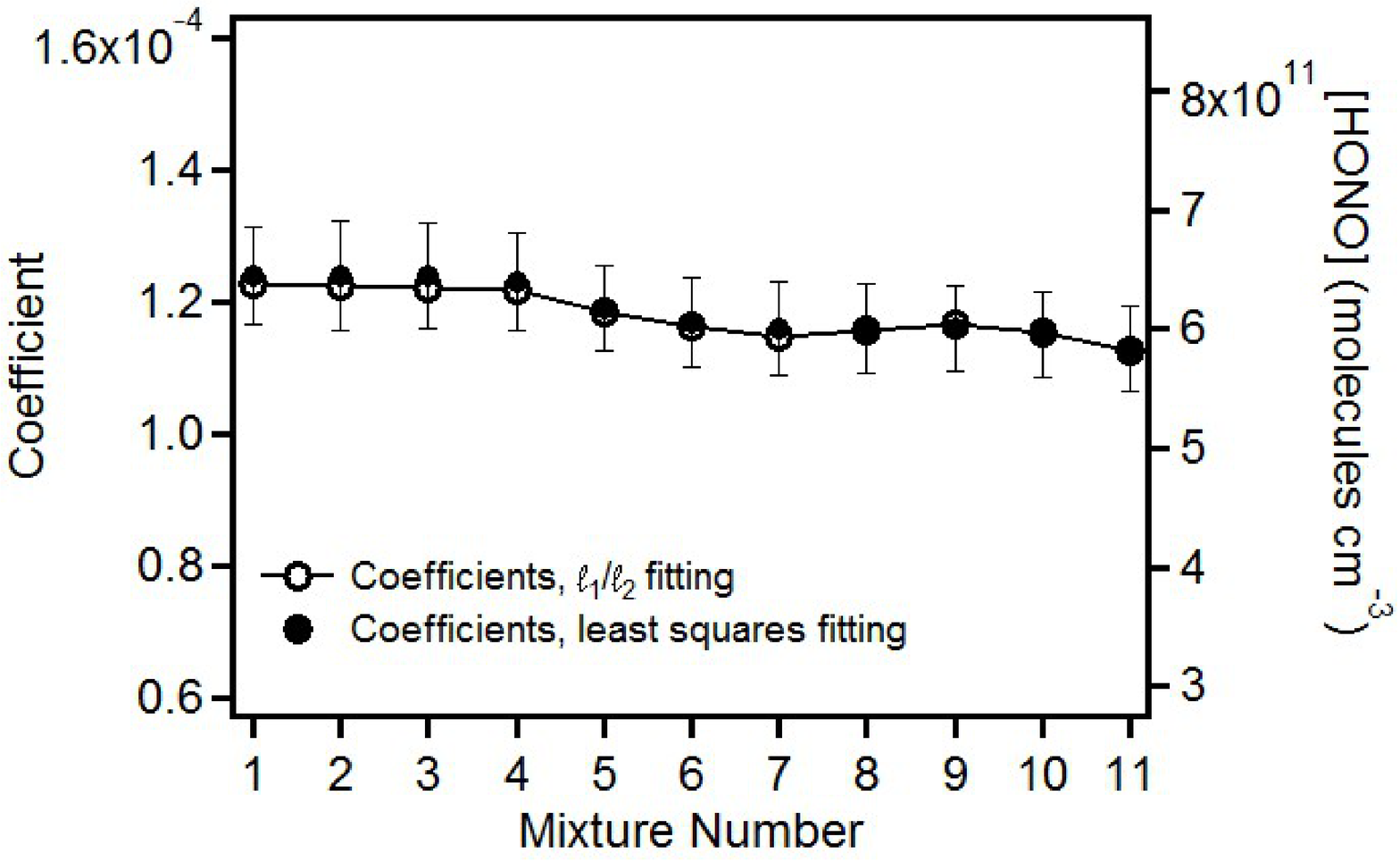}
\hskip 0.5 mm
\includegraphics[height=6cm,width=8cm]{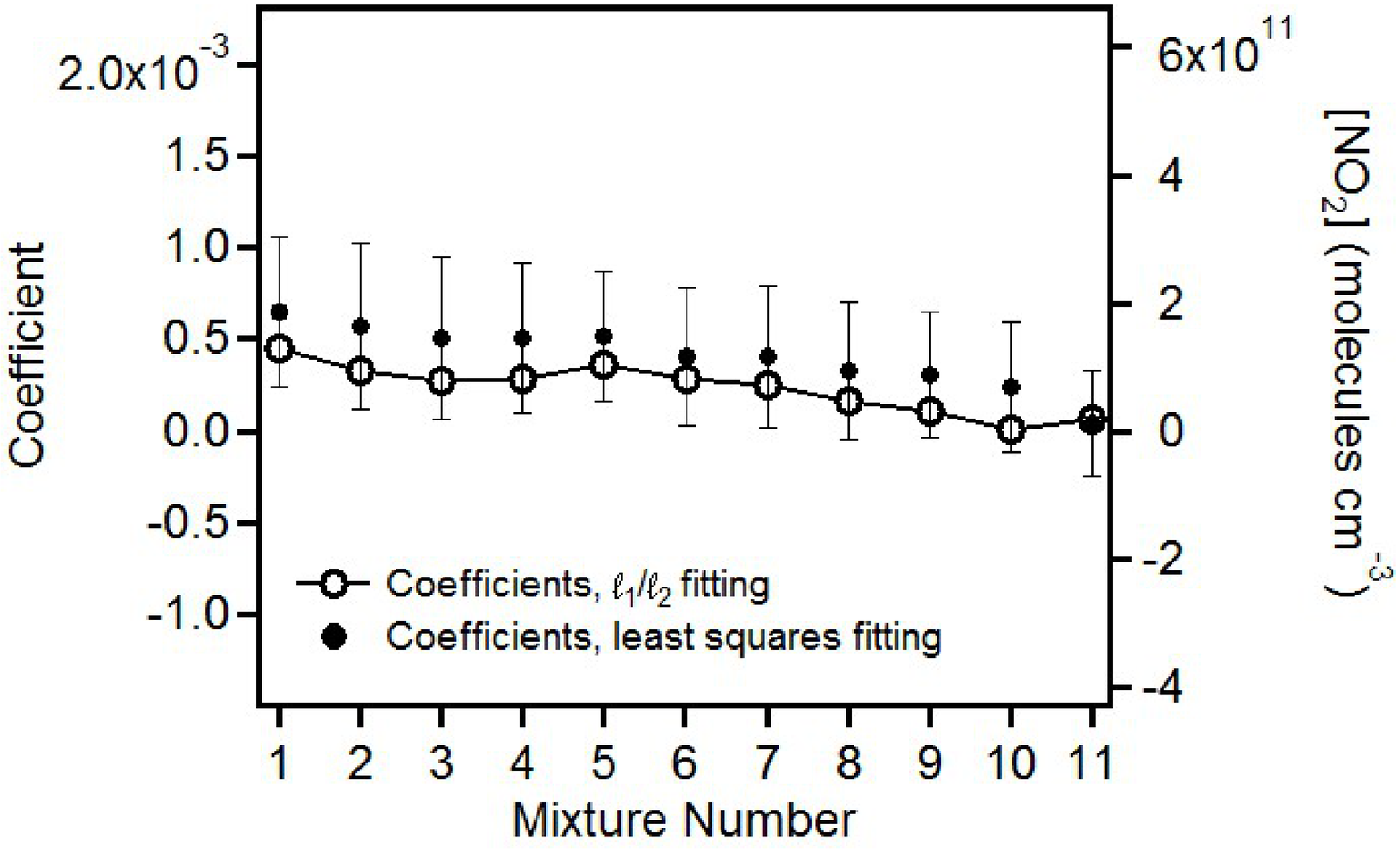}
\includegraphics[height=6cm,width=16cm]{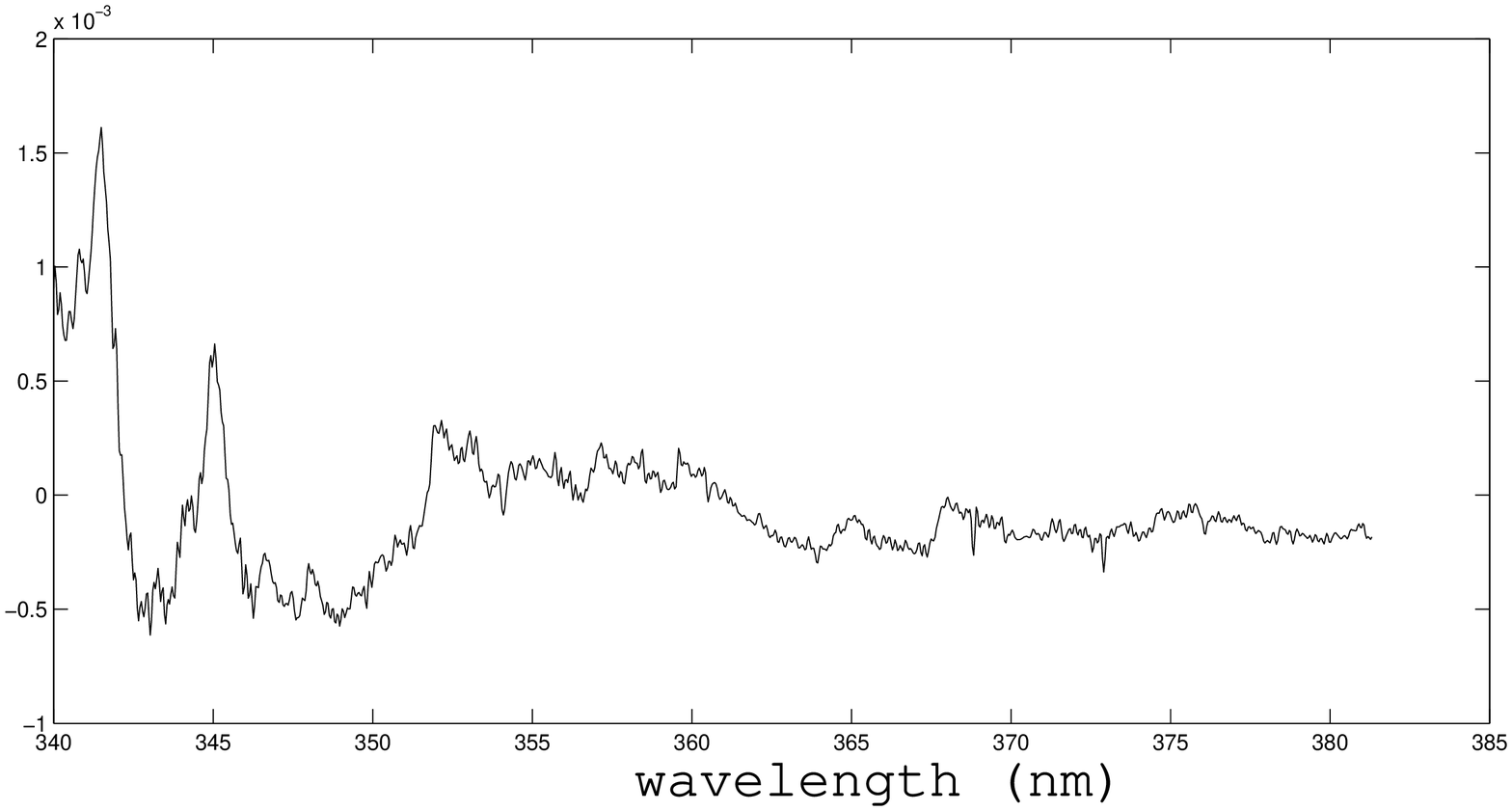}
\caption{(Top row) comparison of the hybrid $\ell_1/\ell_2$ fitting and least squares techniques; HONO coefficients with 2s errors (left), $\mathrm{NO}_2$ coefficients with 2s errors (right), showing good quantitative agreement for eleven mixture spectra collected sequentially over 11 minutes.  Corresponding concentrations in molecules $\mathrm{cm}^{-3}$ are provided on the right axis for each plot. (Bottom row) one fitting residual from robust data fitting. }
 \label{result2}
\end{figure}

\begin{figure}
\includegraphics[height=6.5cm,width=17cm]{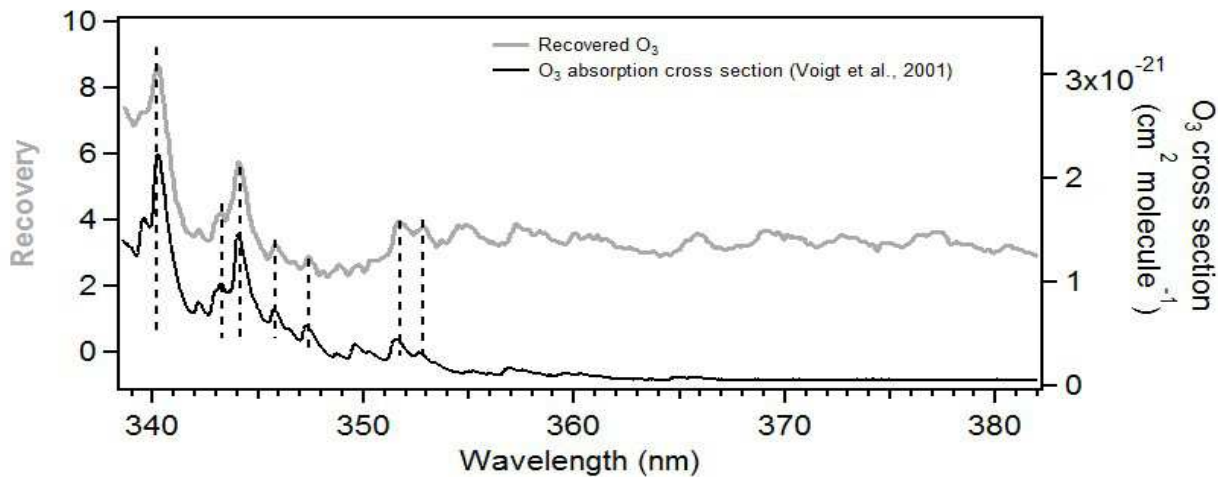}
\caption{Recovered $\mathrm{O}_3$ and its absorption cross section \cite{Voigt_01} for comparison}
 \label{result3}
\end{figure}

We report here the computational results for the proposed method.
In the first example, we fit the known reference spectra of $\mathrm{NO}_2$ and HONO to the DOAS data.
The method identifies an {\em a-priori} unknown trace gas $\mathrm{O}_3$ (ozone) from the fitting residuals.
The results are shown in a series of plots, Fig. \ref{result1}-- Fig. \ref{result3}.
We use 11 sets of data corresponding to different reaction times and hence different gas concentrations  ($X$ has 11 columns) from the experiment. Fig. \ref{result1} illustrates the data preprocessing (EMD) described earlier which removes the fastest and slowest varying components.  The Huber fitting results are presented in Fig. \ref{result2}
which shows the coefficients of $\mathrm{NO}_2$ and HONO in the 11 mixtures in comparison with the concentrations determined using the least squares fitting technique.  These coefficients determined using the hybrid $\ell_1/\ell_2$ fitting technique are all non-negative as well as in very good quantitative agreement with values from least squares fitting.  The fitting residual is in the third plot of Fig. \ref{result2}.
Though some structure can be seen in the residuals, it is not clear if there are other spectral structures embedded in the fitting residuals.  Then further identification was done by JADE.  For the 11 residuals,
we vary the number of independent components in the JADE computation.
We observed that the structure of the first plot in Fig. \ref{result3}
remains approximately invariant as the number varies.
This invariance implies that it should be a hidden trace gas in the fitting residual.
It can be seen that the identified structure resembles $\mathrm{O}_3$ in many peak locations,
especially the region 340--350 nm.  It should be noted that the $\ell_1/\ell_2$ fitting technique currently does not incorporate shifting and squeezing of spectra to optimize fitting, but this can be implemented in the future.  Spectral shifts and squeezes are often used in DOAS analysis routines to account for changes in grating dispersion due to temperature fluctuations and grating positioning accuracy \cite{Platt_Stutz08,Stutz_Platt96}.

\begin{figure}
\includegraphics[height=6cm,width=17cm]{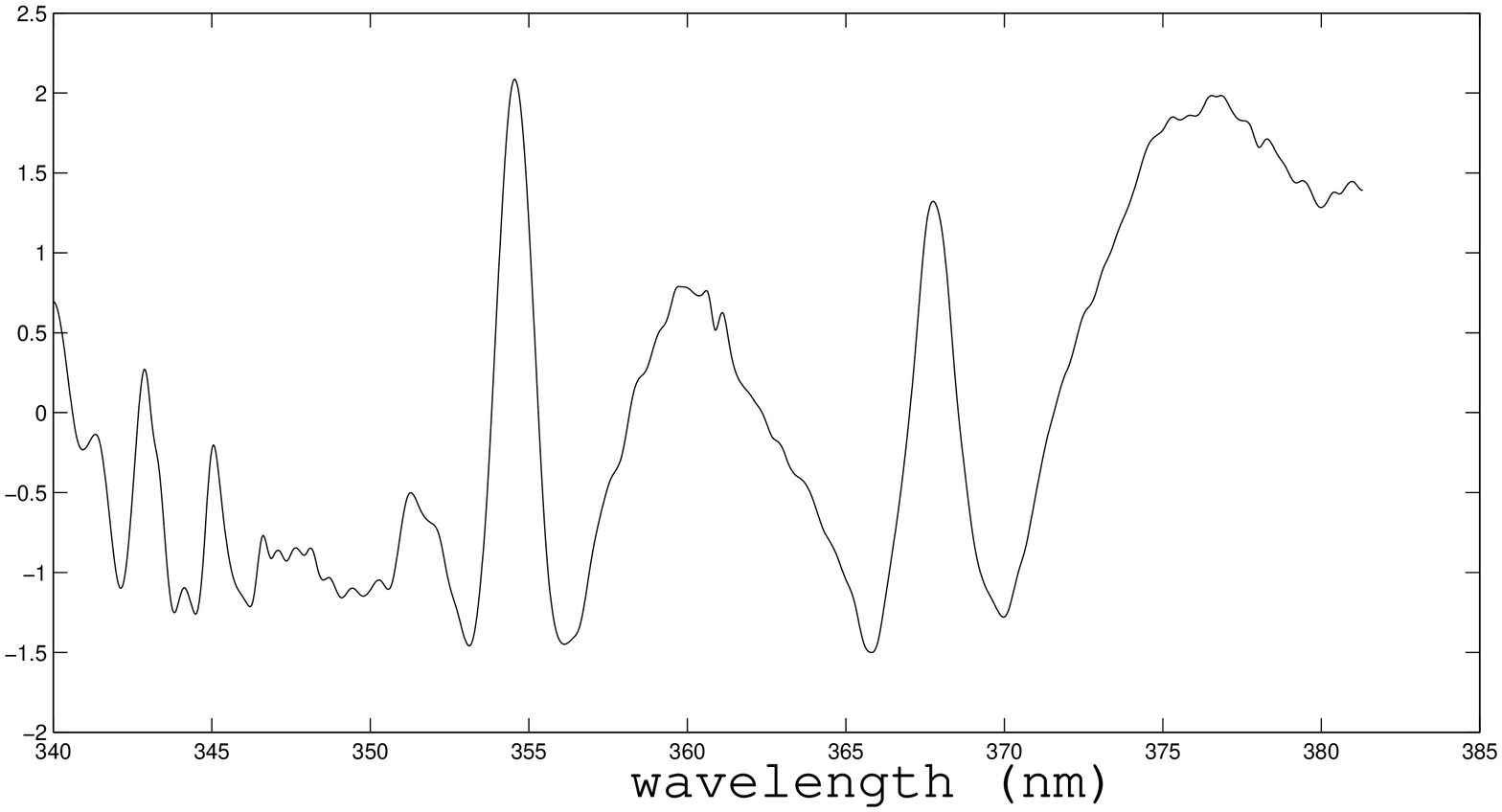}
\includegraphics[height=6cm,width=17cm]{figures/ref_HONO.eps}
\caption{Top plot is the the identified spectral structure 1 compared to the spectral reference of HONO (bottom). }
 \label{result21}
\end{figure}

\begin{figure}
\includegraphics[height=6.5cm,width=17cm]{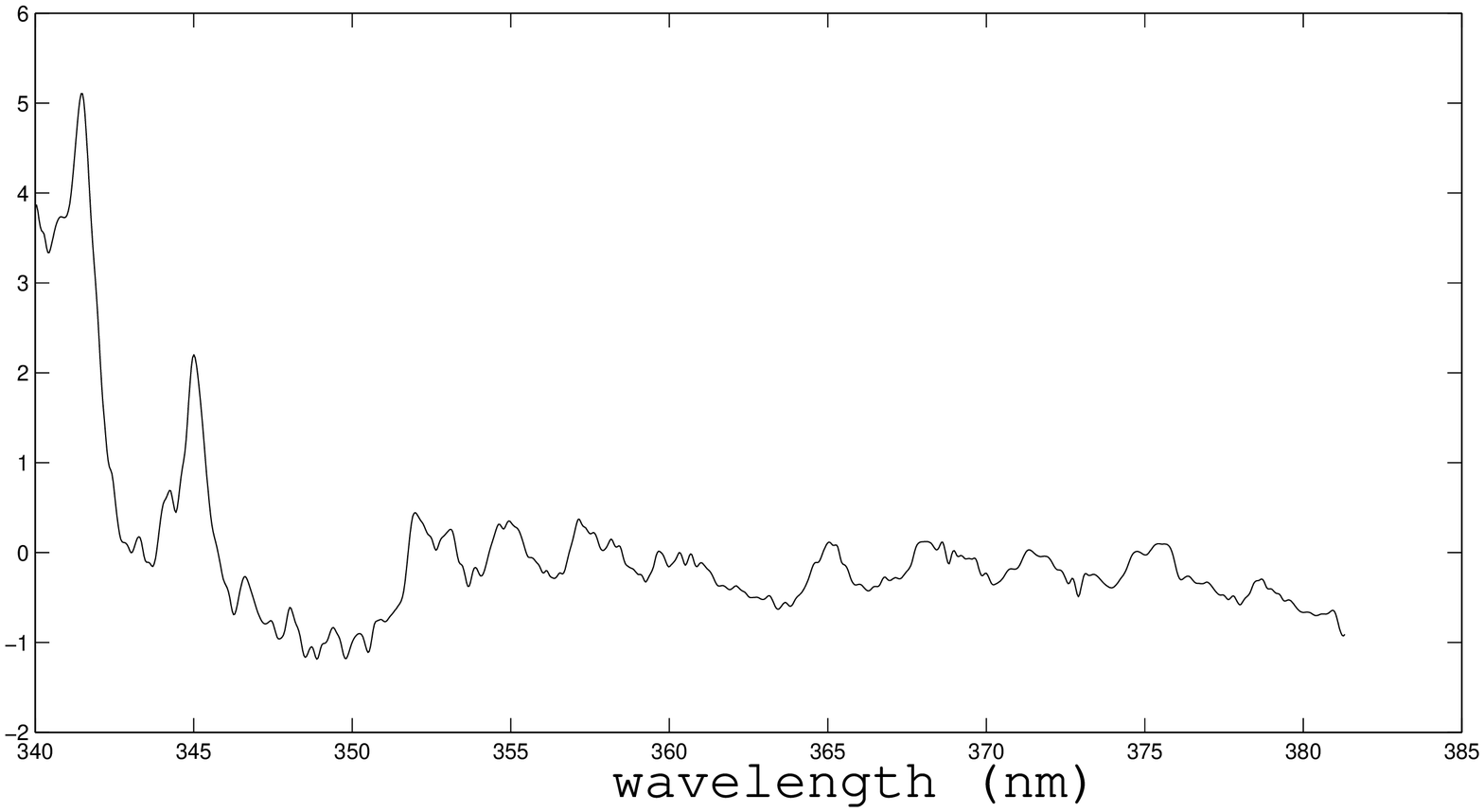}
\includegraphics[height=6.5cm,width=17cm]{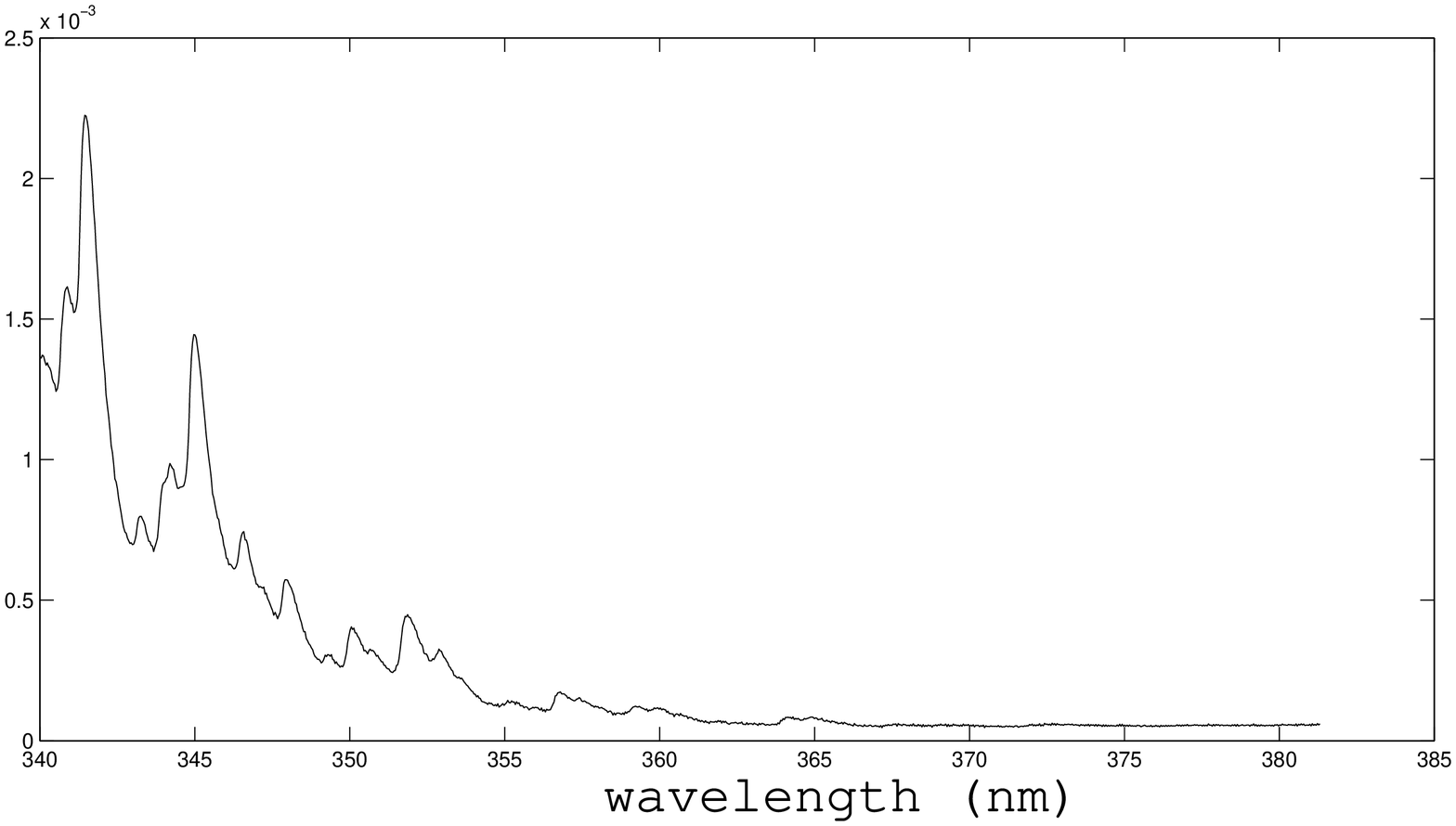}
\caption{Top plot is the the identified spectral structure 2 compared to reference spectrum for $\mathrm{O}_3$ (bottom). }
 \label{result22}
\end{figure}

The second example uses the same set of data (11 mixtures), however
we only use the reference spectrum of $\mathrm{NO}_2$ to fit the data.
Ideally, we should recover $\mathrm{O}_3$ and HONO from the residuals.
The two identified hidden spectral signals are in Fig. \ref{result21} and Fig. \ref{result22}.
The recovered fits are recognizable as HONO and $\mathrm{O}_3$ upon comparison with reference spectra, demonstrating the ability of the technique to identify absorption features {\it without} the use of reference spectra during the fitting procedure. While the $\ell_1/\ell_2$ technique is demonstrated here for laboratory DOAS data with three components, its utility lies in the analysis of atmospheric DOAS spectra, which are more complex.  The least squares method works best when reference spectra for all known absorbing species are used to carry out the fitting, i.e., when the fit residuals are unstructured and do not vary considerably with wavelength.  Given that this condition is rarely satisfied for complex atmospheric measurements, the method described here is complementary in that it can identify species that are either not known to be present or do not yet have available published cross sections.  In addition, the results in Fig. \ref{result2} show its value as an alternative stand-alone technique for analyzing DOAS spectra with the use of appropriate reference spectra.

\section{Concluding Remarks}
We developed a semi-blind source separation method for retrieving the concentrations and performing identifications
of trace gases from DOAS spectra.  The method is designed to identify potentially hidden trace gases
after fitting the known trace gases to the data, which is a challenging problem.
 Our method can be useful for separating unknown source signals
from the residuals after any known reference spectra have been first deployed to fit the data.
 The first novelty of the method is to employ the multi-resolution analysis (EMD) to remove
the slowest varying component from the data.  The removal of such components relies on a polynomial fit in the existing methods. Different polynomials may produce different results, and the degree of the fitting polynomial is often
empirically defined. The multi-resolution approach avoids specifying the order of polynomial,
and it extracts the slow component in an automatic fashion.
The second novelty is to use a hybrid $\ell_1/\ell_2$ interpolated norm (Huber function) to fit the data,
which reduced the effects of outliers and kept the concentrations non-negative.
Lastly, a multi-channel signal decomposition method (JADE) produced encouraging
results on extracting hidden source signals from the fitting residuals.  While use of the least squares fitting procedure for atmospheric data can quantify several trace species simultaneously, typical fit residuals often suggest there are remaining absorbers.  In some cases, species can be inferred based on known atmospheric chemistry, e.g., HONO is often present in $\mathrm{NO}_2$ mixtures.  The major strength of the technique described here is its ability to be used either with existing published reference spectra for quantification or without references for identification of new absorbers.  Numerical results on DOAS data show the promising potential of our method on both trace gas recovery and quantification.

\medskip
This work was partially supported by NSF-ADT grant DMS-0911277 and NSF-CHE 0909227.  The authors thank Professor Jochen Stutz for helpful discussions.

\end{document}